\documentclass[11pt]{scrarticle}

\usepackage[T1]{fontenc}
\usepackage[utf8]{inputenc}
\usepackage{amsfonts,amsmath,amssymb,amsthm,mathtools,dsfont, upgreek,
fullpage,xcolor,hyperref,indentfirst,microtype,authblk}
\usepackage[shortlabels]{enumitem}

\hypersetup{colorlinks,linkcolor={red!50!black},
citecolor={blue!50!black}}

\newcommand{\authorinfo}[3]
{\bigskip
\noindent\textsc{#1:}\\[5pt]#2\\[5pt]
\noindent\textit{Email address: }\href{mailto:#3}{#3}}

\newcommand{\eps}{\varepsilon}
\renewcommand{\alpha}{\upalpha}
\newcommand{\1}{\mathds 1}
\renewcommand{\P}{\mathds P}
\newcommand{\Pc}{\mathcal P}
\newcommand{\E}{\mathds E}
\newcommand{\R}{\mathds R}
\newcommand{\N}{\mathds N}
\newcommand{\Y}{\mathds Y}
\newcommand{\M}{\mathbb M}
\renewcommand{\O}{\mathcal O}
\renewcommand{\o}{{\scriptstyle\mathcal O}}
\newcommand{\B}{\mathcal B}
\newcommand{\T}{\mathcal T}
\newcommand{\G}{\mathcal G}
\renewcommand{\L}{\mathcal L}
\newcommand{\I}{\mathcal I}

\newcommand{\Ec}{\mathcal E}

\newcommand{\Qc}{\mathcal Q}
\newcommand{\hs}{\hspace{1pt}}
\newcommand{\e}{\mathrm e}
\newcommand{\dd}{\mathrm d}
\newcommand{\nto}{\mathrel{\ooalign{$\to$\cr
			\hidewidth\raise.25ex\hbox{\kern-.18em$\scriptscriptstyle/$}\hidewidth\cr}}}
\DeclareMathOperator*{\argmax}{arg\,max}
\DeclareMathOperator{\card}{card}
\DeclareMathOperator{\supp}{supp}
\DeclareMathOperator{\var}{Var}

\DeclareMathOperator{\Leb}{Leb}

\theoremstyle{plain}
\newtheorem{theorem}{Theorem}
\newtheorem{proposition}{Proposition}

\newtheorem{corollary}{Corollary}

\theoremstyle{definition}

\title{Evergreen leaves and universal lifetime laws in preferential-attachment trees}

\author{\textsc{Oleksii Galganov\thanks{Partially supported by a grant from the IMU, supported by the Friends of IMU (FIMU), under the IMU Breakout Fellowship Program.},
Andrii Ilienko}\thanks{Corresponding author.}}

\date{}

\begin{document}
	
\maketitle

\begin{abstract}
	We study leaf lifetimes in general preferential-attachment
	trees. First we characterize the existence of evergreen
	leaves, that is, vertices that remain leaves forever.
	When no evergreen leaves exist, we establish, for a broad
	class of attachment functions, a remarkable one-parameter
	universality: the normalizations and limiting laws for
	various characteristics of leaf lifetimes depend on the
	attachment function only through a single numerical parameter.
	In particular, we obtain a shifted Pareto limit law
	for individual lifetimes and a joint functional limit
	theorem for the lifetime and label of the longest-lasting
	leaf. The limiting two-dimensional process, which we call
	the beta-decorated Fr\'echet extremal process, admits an
	elegant representation in terms of the Stoev--Taqqu
	extremal integral and the sets of optimal choices considered
	by Molchanov and Strokorb.
	We also establish a Poisson point-process limit for rescaled
	leaf lifetimes and corresponding labels, a finiteness criterion
	for long-lasting leaves, and a functional central limit
	theorem for their number at the critical threshold.
\end{abstract}


	
\section{Introduction}

Preferential-attachment models, popularized
by A.-L.~Barab\'asi and R.~Albert \cite{BA99} and subsequently
extended in various directions, have become standard tools
for studying the growth of complex systems, including
the World Wide Web, social networks, and citation graphs.
Their mechanism reflects the well-known ``rich-get-richer''
principle, whereby vertices with more connections are more
likely to attract new ones. Much of the literature has focused
on structural and geometric properties of these networks,
including asymptotic degree distributions \cite{BRST01,RTV07},
maximal degrees \cite{M05}, tree profiles \cite{IK18,KMS17},
as well as local and scaling limits \cite{BBCS14,CDKM15}.
In this paper, we focus on preferential-attachment trees
and study a different, inherently temporal aspect of their
evolution: how long vertices remain leaves and whether
some do so forever.

Consider the general preferential-attachment tree
specified by an arbitrary attachment function
$f:\N_0\to(0,\infty)$.
Starting with the tree $\T_0$ consisting of vertex $0$,
we recursively construct $\T_n$, $n\ge1$, by adding a new vertex $n$
and connecting it to an existing vertex
$k\in\{0,\ldots,n-1\}$, chosen with probability
proportional to $f(\deg_{n-1}^+(k))$:
\begin{equation*}
	\P\{k\to n\mid\T_{n-1}\}
	=\frac{f(\deg_{n-1}^+(k))}
	{\sum_{j=0}^{n-1}f(\deg_{n-1}^+(j))},
	\qquad k=0,\ldots,n-1.
\end{equation*}
Here $\deg_{n-1}^+(k)$ denotes the out-degree of $k$
in $\T_{n-1}$, that is, the number of its children
immediately before vertex $n$ is attached.

The principle ``rich get richer''
also means ``poor stay poor''. Under sufficiently strong
preferential attachment, some vertices may remain leaves
forever; we call such vertices \emph{evergreen leaves}.
It follows from \cite{OS05} that,
for superlinear attachment functions $f(d)\asymp d^p$,
$p>1$, there are a.s.\ infinitely many evergreen leaves.
In Section~\ref{sec:2}, we show that a.s.\ either
no evergreen leaves exist or, in some sense,
almost all vertices remain leaves forever. We give a
necessary and sufficient condition determining which
case occurs, together with explicit, nearly optimal
sufficient conditions on $f$ for each case.

We now turn to the case where no evergreen leaves exist. In this setting, every vertex eventually has a child, but some may remain leaves for exceptionally long periods. This raises several natural questions: how long does a typical vertex
remain a leaf, which vertices do so the longest
and for how long, and how many have unusually long lifetimes?
To address these questions, we impose additional assumptions
on the attachment function $f$. The resulting class remains
very broad, including all sublinear and a large class of
near-linear attachment functions, and covers much of the
regime without evergreen leaves.

It follows from Theorem~2 in \cite{RTV07} that the proportion
of leaves in $\T_n$ converges a.s.\ to $(1+\alpha)^{-1}$,
where $\alpha=f(0)/\lambda^\star$ and $\lambda^\star$ is the
Malthusian parameter of the associated pure-birth
Crump--Mode--Jagers branching process.
Informally, $f(0)$ measures the attractiveness of a leaf
as a potential parent, whereas $\lambda^\star$ represents
the average attachment weight per vertex in the limiting
degree distribution.
In Sections~\ref{sec:3}--\ref{sec:5}, we show that this
dependence on a single parameter extends much further.
Both the normalizations and the limiting laws for
leaf lifetimes, the longest-lasting leaves, and the number
of exceptionally long-lasting leaves are universal up to
the same parameter $\alpha$. In this sense, much of the
fate of leaves is governed by $\alpha$ alone.

We first consider the leaf lifetime $L_k$ of vertex $k$.
In Proposition~\ref{prop:L}, we show that this lifetime is always
of order $k$, and that $L_k/k$ converges in distribution to the
shifted Pareto law with parameter $\alpha$.
Although this follows rather directly from the
branching-process embedding, it is interesting in its own
right and also plays an important role in the subsequent
proofs, allowing us to approximate the dependent leaf lifetimes
by rescaled independent Pareto random variables.

We next study the longest-lasting leaves among the vertices
$k\in\{0,\ldots,n\}$. Set
\begin{equation*}
	I_n=\argmax_{0\le k\le n}L_k,\qquad
	T_n=\max_{0\le k\le n}L_k,
\end{equation*}
that is, $I_n$ stands for the label of a longest-lasting leaf
and $T_n$ for its lifetime. Theorem~\ref{thm:max_argmax} establishes the joint functional convergence in distribution of the processes
$(I_{\lfloor n\cdot\rfloor}/n,
T_{\lfloor n\cdot\rfloor}/n^{1+1/\alpha})$
in the Skorokhod space $D([0,\infty),\R^2)$
to an explicitly described two-dimensional process,
which, in view of its structure, we call the
\emph{beta-decorated $\alpha$-Fr\'echet extremal process}.
Its second component admits an alternative, particularly
elegant representation as an extremal integral with respect
to a Fr\'echet $\sup$-measure in the sense of Stoev and
Taqqu \cite{ST05}, while the first can be described through
the corresponding sets of optimal choices considered by
Molchanov and Strokorb \cite{MS16}.
An important feature of this problem is that the leaf
lifetimes $L_k$ are neither independent nor identically
distributed. Nevertheless, in
Theorem~\ref{thm:point_process_conv}, we establish
the convergence of the point process of rescaled vertex
labels and leaf lifetimes to a Poisson random measure
with explicit intensity.
Besides being of independent interest, this result also provides
the key step in the proof of Theorem~\ref{thm:max_argmax}.
Related questions concerning the $\argmax$ of
non-identically distributed but still independent random
variables were studied in \cite{DR24,DR26}.

The above results show that the maximum leaf lifetime
among the vertices born by time $n$ grows much faster
than the lifetime of a typical vertex. This naturally
leads to the question of how many vertices exhibit
exceptional longevity. In the last section, we study
\emph{$\ell_k$-long-lasting leaves}, that is, vertices $k$
whose leaf lifetimes satisfy $L_k\ge\ell_k$ for a prescribed
sequence of thresholds $(\ell_k)$.
Proposition~\ref{prop:cardL} gives a necessary and sufficient
condition for the number of such vertices to be a.s.\ finite.
In particular, for power thresholds $\ell_k=k^\beta$,
the critical exponent is $\beta=1+1/\alpha$.
At this critical threshold, we show that, on a logarithmic
time scale, the appropriately centered and normalized
counting processes converge in distribution to standard Brownian motion;
see Theorem~\ref{thm:FCLT}.

\section{Evergreen leaves}
\label{sec:2}

We first consider preferential-attachment trees with an arbitrary attachment function $f$. As noted above, we call a vertex $k\in\N_0$ an evergreen leaf if no vertex ever attaches to it as a child. Equivalently, $\deg_n^+(k)=0$ for all $n\ge k+1$. We denote by $\Ec_f$ the (random) set of evergreen leaves.

In the affine preferential-attachment case, that is, when $f(d)=ad+1$ with $a\ge0$, there are a.s.\ no evergreen leaves. Indeed, for any $n\ge k+1$ we have
\begin{equation}
	\label{eq:p_att}
	\P\{k\to n \mid\T_{n-1}\}=\frac{a\deg_{n-1}^+(k)+1}{\sum_{j=0}^{n-1}(a\deg_{n-1}^+(j)+1)}=
	\frac{a\deg_{n-1}^+(k)+1}{a(n-1)+n},
\end{equation}
and hence
\[\P\{k\in\Ec_f\}=\P\left(\bigcap_{n=k+1}^\infty\{k\nto n\}\right)
=\prod_{n=k+1}^\infty\left(1-\frac 1{a(n-1)+n}\right)=0,\qquad k\ge0.\]
Note that the denominator on the right-hand side of \eqref{eq:p_att} is deterministic only because of the affine form of the attachment function.

On the other hand, if $f(d)\asymp d^p$ as $d\to\infty$ for some $p>1$, then $\card\Ec_f=\infty$ a.s. Indeed, by Theorem 1.2 in \cite{OS05} and the remark following it, the limiting tree $\T_\infty=\bigcup_{n=0}^\infty \T_n$ is a.s.\ isomorphic to a tree obtained from a finite tree by attaching to one of its vertices countably many copies of each finite tree with at most $\lfloor (p-1)^{-1}\rfloor+1$ vertices. In particular, $\T_\infty$ contains infinitely many disjoint finite fringe subtrees, each containing a leaf.

Proposition~\ref{prop:EL} below fully characterizes the size of $\Ec_f$. Before stating it, we briefly recall the construction of the pure-birth Crump--Mode--Jagers branching process; see, e.g., our recent paper \cite{GI26}. Consider a pure-birth process with birth rates $f(0),f(1),\ldots$\,, and let $\xi_f$ denote the point process of its jump times. Thus,
\[\xi_f=\sum_{i=0}^\infty\delta_{\tau_0+\ldots+\tau_i},\]
where $\delta_a$ stands for the unit point mass at $a$, and the inter-birth times $\tau_i$ are independent and $\mathsf{Exp}(f(i))$-distributed. 

Let the initial ancestor produce offspring at the atoms of $\xi_f$, and let each individual in turn produce its own offspring according to an independent copy of $\xi_f$ shifted by its birth time. We denote by $\Xi_f$ the pure-birth Crump--Mode--Jagers branching process obtained in this way.
Label the individuals $0,1,\ldots$ in their order of birth, and
denote their birth times by $\sigma_k$, $k\ge0$. The process $\Xi_f$ is said to explode if
\[\sigma_\infty\coloneqq\lim_{k\to\infty} \sigma_k<\infty,\]
that is, if infinitely many individuals are born in finite time. It follows from Proposition 3.5 in \cite{I24} that $\Xi_f$ either explodes a.s.\ or does not explode a.s.

\begin{proposition}~\label{prop:EL}
	\begin{enumerate}[\normalfont(i)]
		\item If $\Xi_f$ does not explode a.s., then
		$\card\Ec_f=0$ a.s.
		\item If $\Xi_f$ explodes a.s., then $\card\Ec_f=\infty$ a.s.
		Even more, almost all vertices are evergreen leaves, in the sense
		that $\Ec_f$ has natural density one{\normalfont:}
		\begin{equation*}
		\lim_{n\to\infty}
		\frac{\card(\Ec_f\cap\{0,\ldots,n\})}{n+1}
		=1\qquad\text{a.s.}
		\end{equation*}
	\end{enumerate}
\end{proposition}

At first sight, the implicit form of Proposition~\ref{prop:EL} may seem to be of limited use, as it reduces one problem to another. But this is not the case. Although explicit necessary and sufficient conditions for explosion of pure-birth Crump--Mode--Jagers branching processes are not known at present, our recent work \cite{GI26} gives close-to-optimal conditions. We now reformulate them in terms of the attachment function $f$ and the set of evergreen leaves $\Ec_f$.

\begin{corollary}~
	\label{cor:EL}
	\begin{enumerate}[\normalfont(i)]
		\item If\, $\sum_{d=0}^\infty \frac 1{f(d)}<\infty$, then $\card\Ec_f=\infty$ a.s.
		\item If
		\begin{equation*}
			\lim_{d\to\infty}\frac{f(d)}d=\infty\qquad\text{and}\qquad
			\liminf_{n\to\infty}\!\!\sum_{d=\lceil\eps\log n\rceil}^n\frac 1{f(d)}>0\quad\text{for some $\eps>0$},
		\end{equation*}
		then $\card\Ec_f=0$ a.s.		
		The same holds if $f$ is bounded above by a function satisfying these two conditions.
		
		In particular, this holds for attachment functions that come very close to the convergence boundary of the series $\sum_{d=0}^\infty \frac 1{f(d)}${\normalfont:}
		\[f(d)=\O\bigl(d\cdot\log d\cdot\log\log d\cdot\ldots\cdot\log^{(l)} d\bigr)\quad\text{for some }l\ge1,\]
		where $\log^{(l)}$ stands for the $l$-fold iterated logarithm.
		\item There exists a function $f$ such that $\sum_{d=0}^\infty \frac 1{f(d)}=\infty$ and $\card\Ec_f=\infty$ a.s.
	\end{enumerate}	
\end{corollary}

\begin{proof}
	This follows from Proposition~\ref{prop:EL} above and Theorem 2.1 in \cite{GI26}.
\end{proof}

Other sufficient conditions for non-explosion of pure-birth Crump--Mode--Jagers branching processes are available in the literature; see, e.g., Remark 2.3 in \cite{GI26}. Via Proposition~\ref{prop:EL}, they also yield sufficient conditions for $\card\Ec_f=0$.

\begin{proof}[Proof of Proposition~\ref{prop:EL}]
	We employ the embedding of the sequence $(\T_n)_{n\ge0}$ into $\Xi_f$, which apparently originates in \cite{P94}; see also the detailed exposition in \cite{OS05}.	
	Denote by $\G_t$ the genealogical tree generated by $\Xi_f$ up to time $t$. Then
	\[(\T_n)_{n\ge0}\overset{\mathrm d}=(\G_{\sigma_n})_{n\ge0},\]
	and hence
	\[\T_\infty\coloneqq\bigcup_{n=0}^\infty\T_n\overset{\mathrm d}=
	\bigcup_{n=0}^\infty \G_{\sigma_n}\eqqcolon\G.\]
	Thus it is enough to prove the claim for $\G$.
	
	If $\Xi_f$ does not explode, then every individual eventually has arbitrarily many children. Hence $\deg_\G^+(k)=\infty$ for all $k\ge0$, and therefore $\card\Ec_f=0$.	Suppose now that $\Xi_f$ explodes. In this case, the claim is
	implicit in the proof of Theorem~6.11 in \cite{HJ17}.
	Although that theorem additionally assumes $\sum_{d=0}^\infty1/f(d)<\infty$, its proof actually requires only $\sigma_\infty<\infty$. As the argument can be given in just a few lines, we include it
	for completeness.
	
	Let $\delta>0$. Since the first reproduction times $E_k$, measured
	from the individuals' births, are independent and $\mathsf{Exp}(f(0))$-distributed, the strong law of large numbers yields
	\[\frac{\sum_{k=0}^n\1\{E_k>\delta\}}{n+1}\to\e^{-\delta f(0)}
	\qquad\text{a.s.}\]
	All but finitely many
	individuals are born after $\sigma_\infty-\delta$, and those with
	$E_k>\delta$ remain childless up to explosion. Therefore,
	\[\liminf_{n\to\infty}
	\frac{\card(\Ec_f\cap\{0,\ldots,n\})}{n+1}
	\ge\e^{-\delta f(0)}\qquad\text{a.s.}\]
	Since $\delta>0$ is arbitrary, the claim follows.
\end{proof}

\section{Longest-lasting leaves}
\label{sec:3}

In this section, we deal with the non-explosive regime for the associated
branching process $\Xi_f$. By
Proposition~\ref{prop:EL}, the tree
$(\T_n)_{n\ge0}$ then has no evergreen leaves. This leads to the following
natural questions: which vertices remain leaves the longest, and for
how long?

We impose two standard conditions on the attachment function $f$
and start with the Malthusian condition \eqref{eq:M}; see, e.g., condition~(M)
in \cite{RTV07}. Let
\[
\hat f(\lambda)
=\int_0^\infty \e^{-\lambda t}\,(\E\xi_f)(\dd t)
=\sum_{k=1}^\infty
\prod_{i=0}^{k-1}\frac{f(i)}{f(i)+\lambda}
\]
be the Laplace transform of the mean measure $\E\xi_f$ of the pure-birth
point process $\xi_f$, and put
\[
\underline\lambda
=\inf\{\lambda>0:\hat f(\lambda)<\infty\}.
\]
We assume that
\begin{equation}
	\label{eq:M}
	\underline\lambda<\infty\quad\text{and}\quad
	\lim_{\lambda\hs\downarrow\hs\underline\lambda}
	\hat f(\lambda)>1.
\end{equation}
Under \eqref{eq:M}, there exists a unique
$\lambda^\star>\underline\lambda$ such that
$\hat f(\lambda^\star)=1$. It is called the \emph{Malthusian parameter}
of $\Xi_f$. For the affine attachment function $f(d)=ad+1$, one has
$\lambda^\star=a+1$; see, e.g., Section~4.2 in \cite{RTV07}. In general,
$\lambda^\star$ can be computed numerically.

We also impose the well-known $x\log x$ condition; see, e.g.,
Section~5 in \cite{JN84}. Set
\[
\hat\xi_f(\lambda^\star)
=\int_0^\infty \e^{-\lambda^\star t}\,\xi_f(\dd t)
\]
and assume that
\begin{equation}
	\label{eq:L}
	\E\bigl(\hat\xi_f(\lambda^\star)
	\log^+\hat\xi_f(\lambda^\star)\bigr)<\infty.
\end{equation}
Condition \eqref{eq:L} is necessary and sufficient for the limiting random
variable in the Malthusian growth theorem to be a.s.\ positive; see
Theorem~5.3 in \cite{JN84}.

In what follows, \eqref{eq:M} and \eqref{eq:L} will be our standing
assumptions. They hold for a broad class of attachment functions,
covering much of the non-explosive regime. Several explicit sufficient
conditions are available in the literature. The Malthusian condition
\eqref{eq:M} holds in the affine case $f(d)=ad+1$ and for strongly
sublinear attachment functions; see Section~4.2 in \cite{RTV07} and
Example~4 in \cite{GVCH17}. Under \eqref{eq:M}, condition \eqref{eq:L}
is satisfied whenever $f(d)\to\infty$; see Lemma~1 in \cite{RTV07}. We now give simple explicit conditions, hopefully also of independent interest, ensuring
\eqref{eq:M}$+$\eqref{eq:L} and extending all these results.

\begin{proposition}
	\label{prop:cond}
	Suppose that one of the following conditions holds as $d\to\infty${\normalfont:}
	\begin{enumerate}[\normalfont(i)]
		\item $f(d)=ad+\O(d^\gamma)$ for some $a>0$ and $\gamma<1$
		(the near-linear case);
		\item $f(d)=\o(d)$
		(the sublinear case).
	\end{enumerate}
	Then both \eqref{eq:M} and \eqref{eq:L} are satisfied.
\end{proposition}

\begin{proof}
	Assume first that \textup{(i)} holds. Since $\log(1+x)=x+\O(x^2)$ and
	$\sum_{i=0}^\infty f(i)^{-2}<\infty$,
	\[\log\prod_{i=0}^{k-1}\frac{f(i)}{f(i)+\lambda}
	=-\sum_{i=0}^{k-1}\log\biggl(1+\frac{\lambda}{f(i)}\biggr)
	=-\lambda\sum_{i=0}^{k-1}\frac1{f(i)}+\O(1).\]
	Since $f(i)^{-1}-(ai)^{-1}=\O(i^{\gamma-2})$, this in turn equals
	$-\frac{\lambda}{a}\log k+\O(1)$.
	Hence,
	\[\prod_{i=0}^{k-1}\frac{f(i)}{f(i)+\lambda}
	\asymp k^{-\lambda/a},\qquad k\to\infty.\]
	Therefore, $\hat f(\lambda)<\infty$ if and only if
	$\lambda>a$. So, $\underline\lambda=a$ and
	$\lim_{\lambda\downarrow a}\hat f(\lambda)=\infty$ by monotone convergence. Thus, \eqref{eq:M} holds.
	
	In case \textup{(ii)},
	$f(i)\le\lambda i/2$ for large $i$, and hence
	$f(i)/(f(i)+\lambda)\le i/(i+2)$. Hence,
	\[\prod_{i=0}^{k-1}\frac{f(i)}{f(i)+\lambda}
	=\O\Biggl(\prod_{i=1}^{k-1}\frac i{i+2}\Biggr)
	=\O(k^{-2}).\]
	Thus, $\hat f(\lambda)<\infty$ for every $\lambda>0$, and so
	$\underline\lambda=0$. Moreover,
	$\lim_{\lambda\downarrow0}\hat f(\lambda)=\infty$ by monotone convergence.
	Thus, \eqref{eq:M} holds again.
	
	We now verify \eqref{eq:L}. In case \textup{(i)},
	$\lambda^\star>\underline\lambda=a$, whereas in case \textup{(ii)},
	$\lambda^\star>0$. Thus, in either case, there exist
	$c\in(0,\lambda^\star)$ and $b>0$ such that
	$f(i)\le ci+b$ for all $i\ge0$. Let $\tau_i\sim\mathsf{Exp}(f(i))$ be the independent	inter-birth times of $\xi_f$. Writing $\|X\|_2$ for $\sqrt{\E X^2}$, we have by Minkowski's inequality
	\begin{align*}
		\bigl\|\hat\xi_f(\lambda^\star)\bigr\|_2
		&=\Bigl\|\sum_{k=1}^\infty
		\e^{-\lambda^\star(\tau_0+\ldots+\tau_{k-1})}\Bigr\|_2
		\le\sum_{k=1}^\infty
		\bigl\|\e^{-\lambda^\star(\tau_0+\ldots+\tau_{k-1})}\bigr\|_2\\
		&=\sum_{k=1}^\infty
		\prod_{i=0}^{k-1}
		\bigl(\E\e^{-2\lambda^\star\tau_i}\bigr)^{1/2}
		=\sum_{k=1}^\infty
		\prod_{i=0}^{k-1}
		\biggl(\frac{f(i)}{f(i)+2\lambda^\star}\biggr)^{1/2}
		\le\sum_{k=1}^\infty
		\prod_{i=0}^{k-1}
		\biggl(\frac{ci+b}{ci+b+2\lambda^\star}\biggr)^{1/2}\!.
	\end{align*}
	By the standard asymptotics for ratios of Gamma functions, the last
	product is $\O(k^{-\lambda^\star/c})$. Since
	$c<\lambda^\star$, the series converges. Thus,
	$\E\hat\xi_f(\lambda^\star)^2<\infty$, which implies \eqref{eq:L}.
\end{proof}

Before turning to the longest-lasting leaves, we first look at the lifetime
$L_n$ of a single vertex $n$ as a leaf. Surprisingly, throughout the entire
class of preferential-attachment trees satisfying \eqref{eq:M} and
\eqref{eq:L}, this lifetime is always of order $n$. Even more, after
normalization by $n$, the limit law is universal up to a single numerical
parameter determined by the attachment function $f$.

\begin{proposition}
	\label{prop:L}
	Let $P_\alpha$ have the shifted Pareto distribution with tail
	\begin{equation}
	\label{eq:Pareto}
	\P\{P_\alpha>x\}=(1+x)^{-\alpha},\qquad x\ge0.
	\end{equation}
	Then, under \eqref{eq:M} and \eqref{eq:L},
	\[\frac{L_n}{n}\xrightarrow{\mathrm d} P_\alpha,\qquad n\to\infty,\]
	where $\alpha=f(0)/\lambda^\star$.
\end{proposition}

Besides being of independent interest, Proposition~\ref{prop:L} will also
be instrumental in the subsequent proofs, where rescaled independent copies of $P_\alpha$ will serve as approximations to the lifetimes $L_n$.

\begin{proof}[Proof of Proposition \ref{prop:L}]
	Denote by $\sigma_n$ the birth time of individual $n$ in $\Xi_f$ and by
	$N(t)$ its population size at time $t$. By the standard Malthusian growth theorem, say, in the form of	Proposition~3 in \cite{GVCH17},
	\begin{equation}
	\label{eq:Mgr}
	\e^{-\lambda^\star t}N(t)=W+\o(1)\quad\text{a.s.},
	\qquad t\to\infty,
	\end{equation}
	for some non-negative random variable $W$; the assumptions of that
	proposition are readily verified under \eqref{eq:M}. As noted above,
	condition \eqref{eq:L} additionally implies $W>0$ a.s.
	
	Since $N(\sigma_n)=n+1$, we obtain
	\begin{equation}
	\label{eq:siglog}
	\lambda^\star\sigma_n=\log n-\log W+\o(1)\quad\text{a.s.},
	\qquad n\to\infty.
	\end{equation}
	Let $E_n\sim\mathsf{Exp}(f(0))$ be the first reproduction time of
	individual $n$, measured from its
	birth. Then
	\[\lambda^\star E_n=\lambda^\star\sigma_{n+L_n}-\lambda^\star\sigma_n
	=\log(1+L_n/n)+\o(1)\quad\text{a.s.}\]
	Hence, $L_n/n\xrightarrow{\mathrm d}\e^{\lambda^\star E}-1$,
	where $E\sim\mathsf{Exp}(f(0))$, and a straightforward calculation shows
	that the right-hand side has the shifted Pareto distribution
	\eqref{eq:Pareto} with $\alpha=f(0)/\lambda^\star$.
\end{proof}

We now turn to the longest-lasting leaves up to time $\lfloor nt\rfloor$ and establish a
joint functional limit theorem for the label of such a leaf and its lifetime.
As we shall see, the maximal lifetime may grow at a substantially faster rate
than $n$, the order appearing in Proposition~\ref{prop:L}. Remarkably, both
this rate and the limiting two-dimensional process depend on the attachment
function $f$ once again only through the same parameter
$\alpha=f(0)/\lambda^\star$.

We first define the two-dimensional process $\bigl(X_\alpha(t),Y_\alpha(t)\bigr)_{t\ge 0}$, which will appear as the limit. In view of its structure, we call it the beta-decorated $\alpha$-Fr\'echet extremal process. The second component $Y_\alpha$ is a modified Fr\'echet extremal process; for general background, see, e.g., Section~4.3 in \cite{R08}. The process $Y_\alpha$ is non-decreasing, pure-jump and c\`adl\`ag, with $Y_\alpha(0)=0$, and is characterized by the following recursion: for $0\le s<t$,
\begin{equation}
	\label{eq:Y_rec}
	Y_\alpha(t)=\max\bigl\{Y_\alpha(s),\|g_{s,t}\|_\alpha\hs F_{\alpha,s,t}\bigr\},
\end{equation}
where
\[g_{s,t}(u)=u\1_{(s,t]}(u),\qquad
\|g_{s,t}\|_\alpha=\left(\int_0^\infty g_{s,t}(u)^\alpha\,\dd u\right)^{1/\alpha}
=\left(\frac{t^{\alpha+1}-s^{\alpha+1}}{\alpha+1}\right)^{1/\alpha}\!,\]
and, for each $0\le s<t$, $F_{\alpha,s,t}$ is a standard $\alpha$-Fr\'echet random variable, independent of $Y_\alpha(s)$, with distribution function
\[\P\{F_{\alpha,s,t}\le x\}=\e^{-x^{-\alpha}},\qquad x>0.\]
Moreover, the variables $F_{\alpha,s,t}$ corresponding to pairwise disjoint intervals $(s,t]$ are independent.

Let $d\ge1$, $0=t_0<t_1<\ldots<t_d$ and set $X_\alpha(0)=0$. The conditional distribution of $\bigl(X_\alpha(t_1),\ldots,X_\alpha(t_d)\bigr)$ given $\bigl(Y_\alpha(t_1),\ldots,Y_\alpha(t_d)\bigr)$ is defined recursively by
\begin{equation}
\label{eq:X_dist}
X_\alpha(t_i)=
\begin{cases}
	X_\alpha(t_{i-1}),&\text{if } Y_\alpha(t_i)=Y_\alpha(t_{i-1}),\\
	B_{\alpha,t_{i-1},t_i},&\text{if } Y_\alpha(t_i)>Y_\alpha(t_{i-1}),
\end{cases}
\qquad i=1,\ldots,d,
\end{equation}
where $B_{\alpha,t_{i-1},t_i}$ are independent of each other and of $(Y_\alpha(t_1),\ldots,Y_\alpha(t_d))$ and have
beta-type densities
\[f_{B_{\alpha,s,t}}(x)=\|g_{s,t}\|_\alpha^{-\alpha}x^\alpha\1_{(s,t]}(x),
\qquad 0\le s<t.\]
The recursions \eqref{eq:Y_rec} and \eqref{eq:X_dist} uniquely determine the finite-dimensional distributions of $(X_\alpha,Y_\alpha)$.

The process $(X_\alpha,Y_\alpha)$ admits a more elegant description in terms of the extremal stochastic integrals ${\vphantom\int}^e\hspace{-5pt}\int$ of Stoev and Taqqu \cite{ST05} as well as the corresponding sets of optimal
choices considered by Molchanov and Strokorb \cite{MS16}. Let $\Leb$ denote Lebesgue measure, and, for $\alpha>0$, let $\M_{\alpha}$ be a standard $\alpha$-Fr\'echet $\sup$-measure on $\B([0,\infty))$, that is, a $\sigma$-$\sup$-additive random set function with independent values on disjoint sets such that, for every $B\in\B([0,\infty))$ with $\Leb(B)<\infty$,
\[\P\{\M_{\alpha}(B)\le x\}=\e^{-\Leb(B)\hs x^{-\alpha}},\qquad x>0;\]
see Definition~2.1 in \cite{ST05}. Then $Y_\alpha$ can be realized as
\[Y_\alpha(t)={\vphantom\int}^e\hspace{-8pt}\int_0^t u\,\M_\alpha(\dd u),\qquad t\ge0.\]
Indeed, recursion \eqref{eq:Y_rec} follows from the properties of the Stoev--Taqqu integral: for $0\le s<t$,
\[Y_\alpha(t)=\max\biggl\{
{\vphantom\int}^e\hspace{-8pt}\int_0^\infty u\1_{(0,s]}(u)\,\M_\alpha(\dd u),
{\vphantom\int}^e\hspace{-8pt}\int_0^\infty u\1_{(s,t]}(u)\,\M_\alpha(\dd u)
\biggr\}.\]
The two terms in the maximum are independent, the first one equals $Y_\alpha(s)$, and the second one has the distribution of $\|g_{s,t}\|_\alpha F_\alpha$, where $F_\alpha$ is a standard $\alpha$-Fr\'echet; see \cite[p.~250]{ST05}.

Now, for $0\le s<t$, let
\[\Y_{\alpha,s,t}(B)={\vphantom\int}^e\hspace{-8pt}\int_B u\,\M_\alpha(\dd u),\qquad B\in\B((s,t]).
\]
Then $\Y_{\alpha,s,t}$ is a completely random max-stable
$\sup$-measure on $(s,t]$. Hence, by Theorem 8.2 in \cite{MS16}, the set of optimal choices
\[\bigl\{u\in(s,t]:
\Y_{\alpha,s,t}(\{u\})=\Y_{\alpha,s,t}((s,t])\bigr\}\]
is a.s.~a singleton, whose element is independent of
$\Y_{\alpha,s,t}((s,t])$. This element, which has the required distribution, is taken as
$B_{\alpha,s,t}$ in the above construction.

Finally, we note a third representation of $(X_\alpha,Y_\alpha)$
established and used in the proof. Let $\Lambda$ be a Poisson
random measure on $[0,\infty)\times(0,\infty)$ with intensity
$\alpha x^\alpha y^{-\alpha-1}\,\dd x\,\dd y$.
Then $(X_\alpha(t),Y_\alpha(t))$ can be realized as
the a.s.\ unique highest atom of $\Lambda$ in
$(0,t]\times(0,\infty)$.

The processes $X_\alpha$ and $Y_\alpha$ are not independent,
although $X_\alpha(t)$ and $Y_\alpha(t)$ are independent for
each fixed $t$. The latter follows either from \eqref{eq:X_dist}
by taking $d=1$ and $t_1=t$, or from Theorem~8.2 in \cite{MS16}.
In particular,
\begin{equation}
	\label{eq:BF}
	\bigl(X_\alpha(1),Y_\alpha(1)\bigr)
	\overset{\mathrm d}{=}
	\bigl(B_\alpha,(\alpha+1)^{-1/\alpha}F_\alpha\bigr),
\end{equation}
where $B_\alpha\sim\mathsf{Beta}(\alpha+1,1)$ and
$F_\alpha\sim\mathsf{Frechet}(\alpha)$ are independent.

\begin{theorem}\label{thm:max_argmax}
	Let
	\[I_n=\argmax_{0\le k\le n}L_k,\qquad
	T_n=\max_{0\le k\le n}L_k,\]
	where $\argmax$ stands for the least maximizer. Then, under \eqref{eq:M} and \eqref{eq:L},
	\[\left(\frac{I_{\lfloor n\cdot\rfloor}}n,
	\frac{T_{\lfloor n\cdot\rfloor}}{n^{1+1/\alpha}}\right)
	\overset{\mathrm d}\to(X_\alpha,Y_\alpha),\qquad n\to\infty,\]
	in the $J_1$ topology on $D\bigl([0,\infty),\R^2\bigr)$, where $\alpha=f(0)/\lambda^\star$ and $(X_\alpha,Y_\alpha)$ is the beta-decorated $\alpha$-Fr\'echet extremal process.
\end{theorem}

\begin{corollary}
	\label{cor:IT}
	$I_n/n$ and $T_n/n^{1+1/\alpha}$ are asymptotically independent and converge in distribution to $B_\alpha$ and $(\alpha+1)^{-1/\alpha}F_\alpha$ from \eqref{eq:BF}, respectively.
\end{corollary}

The difference between the normalizations of $T_n$ in
Corollary~\ref{cor:IT} and $L_n$ in Proposition~\ref{prop:L} can be
explained by the fact that only a vanishing fraction of the lifetimes $L_k$, $k=1,\ldots,n$, exceed any superlinear level. We illustrate this
in the affine case $f(d)=ad+1$, $a\ge0$, where
$\alpha=f(0)/\lambda^\star=(a+1)^{-1}$,
and all calculations can be carried out explicitly.

For $k,l\ge1$, by \eqref{eq:p_att},
\begin{align*}
	\P\{L_k>l\}&=\P\{k\nto k+1,\ldots,k\nto k+l\}\\
	&=\prod_{j=k+1}^{k+l}\biggl(1-\frac{1}{a(j-1)+j}\biggr)
	=\prod_{j=k+1}^{k+l}\frac{j-1}{j-\frac{a}{a+1}}
	=\frac{\Gamma(k+\alpha)}{\Gamma(k)}
	\frac{\Gamma(k+l)}{\Gamma(k+l+\alpha)}.
\end{align*}
Therefore, for every $\varepsilon>0$,
\[\E\sum_{k=1}^n\1\{L_k>n^{1+\eps}\}=
\sum_{k=1}^n\frac{\Gamma(k+\alpha)}{\Gamma(k)}
\frac{\Gamma\bigl(k+\lfloor n^{1+\eps}\rfloor\bigr)}
{\Gamma\bigl(k+\lfloor n^{1+\eps}\rfloor+\alpha\bigr)}
\sim n^{-\alpha(1+\eps)}
\sum_{k=1}^n\frac{\Gamma(k+\alpha)}{\Gamma(k)},\qquad n\to\infty.\]
The last sum is $\sim\frac{1}{1+\alpha}n^{1+\alpha}$
by the Stolz--Ces\`aro theorem. Hence, the expected number of exceedances is
$\sim\frac{1}{1+\alpha}n^{1-\alpha\eps}=\o(n)$.
Moreover, it is of constant order when $\eps=\alpha^{-1}$, that is, at
the scale $n^{1+1/\alpha}$,
which explains the normalization of $T_n$. A more detailed analysis of
the number of exceedances is given in
Section~\ref{sec:5}.

\section{Proof of Theorem~\ref{thm:max_argmax}}

The proof is based on the following result of independent interest.
Let $Q=[0,\infty)\times(0,\infty)$ and equip it with the localizing ring
\begin{equation}
\label{eq:Q}
\Qc=\bigl\{B\in\B(Q): B\subset [0,T]\times[\eps,\infty)
\text{ for some } T,\eps>0\bigr\};
\end{equation}
see, e.g., \cite[p.~19]{K17} for background on localizing rings. On this
space, consider the random point measure
\begin{equation}
	\label{eq:Lambda_n}
	\Lambda_n=\sum_{k=0}^\infty
	\delta_{\left(k/n,L_k/n^{1+1/\alpha}\right)}.
\end{equation}
Recall that the vague topology on the space of locally finite measures on
$(Q,\Qc)$ is generated by the integration maps
\[
\nu\mapsto\int_Q g\,\dd\nu,
\]
where $g$ ranges over all bounded continuous functions whose support belongs to
$\Qc$; see, e.g., Section~3.4 in \cite{R08} or Chapter~4 in \cite{K17}
for a general exposition. Following \cite{K17}, we denote vague convergence by $\xrightarrow{\mathrm v}$ and convergence in distribution w.r.t.~the vague topology by
$\xrightarrow{\mathrm{vd}}$.

\begin{theorem}
	\label{thm:point_process_conv}
	Let $\Lambda$ be a
	Poisson random measure with intensity
	\[\mu(\dd x,\dd y)=\alpha x^\alpha y^{-\alpha-1}\,\dd x\,\dd y.\]
	Then, under \eqref{eq:M} and \eqref{eq:L},
	$\Lambda_n\xrightarrow{\mathrm{vd}}\Lambda$.
\end{theorem}

\begin{proof}
	The idea of the proof is to approximate $\Lambda_n$ by a simpler point
	measure $\widetilde\Lambda_n$ with independent atoms generated by
	i.i.d.\ copies of $P_\alpha$ from
	Proposition~\ref{prop:L}, prove the convergence for
	$\widetilde\Lambda_n$, and then transfer it to $\Lambda_n$.
	
	Denote these copies by $P_{\alpha,k}$, $k\ge0$. By analogy
	with \eqref{eq:Lambda_n}, define the random point measure
	\[\widetilde\Lambda_n=\sum_{k=0}^\infty
	\delta_{\left(k/n,(k+1)P_{\alpha,k}/n^{1+1/\alpha}\right)}
	\eqqcolon\sum_{k=0}^\infty\eta_{\alpha,k,n}.\]
	We first show that
	$\widetilde\Lambda_n\xrightarrow{\mathrm{vd}}\Lambda$.
	For $0\le s<t$ and $y>0$, as $n\to\infty$,
	\begin{multline*}
		\E\widetilde\Lambda_n\bigl((s,t]\times(y,\infty)\bigr)
		=\sum_{sn<k\le tn}
		\P\biggl\{\frac{(k+1)P_{\alpha,k}}{n^{1+1/\alpha}}>y\biggr\}\\
		=\sum_{sn<k\le tn}
		\biggl(1+\frac{y n^{1+1/\alpha}}{k+1}\biggr)^{-\alpha}
		\sim\frac{y^{-\alpha}}{n}
		\sum_{sn<k\le tn}
		\left(\frac{k+1}{n}\right)^\alpha
		\to
		y^{-\alpha}\int_s^t x^\alpha\,\dd x
		=\mu\bigl((s,t]\times(y,\infty)\bigr),
	\end{multline*}
	which implies
	$\E\widetilde\Lambda_n\xrightarrow{\mathrm v}\mu=\E\Lambda$.
	Furthermore, for every $B\in\Qc$, as $n\to\infty$,
	\[\sup_{k\ge0}\P\{\eta_{\alpha,k,n}(B)>0\}\le
	\max_{0\le k\le Tn}\P\biggl\{\frac{(k+1)P_{\alpha,k}}{n^{1+1/\alpha}}\ge\eps
	\biggr\}=\max_{0\le k\le Tn}
	\biggl(1+\frac{\eps n^{1+1/\alpha}}{k+1}\biggr)^{-\alpha}\to0,\]
	where $\eps$ and $T$ are as in \eqref{eq:Q}. Thus, since
	$\eta_{\alpha,k,n}$ are independent, they form a null array in the sense
	of Section~4.3 in \cite{K17}. Hence, the Grigelionis theorem in the form
	of Corollary~4.25 in \cite{K17} yields	$\widetilde\Lambda_n\xrightarrow{\mathrm{vd}}\Lambda$.
	
	It remains to transfer this convergence to $\Lambda_n$. As in the proof
	of Proposition~\ref{prop:L}, let $\sigma_k$ be the birth time of individual
	$k$ in $\Xi_f$ and $E_k$ its first reproduction time, measured from its
	birth. Thus, we may take
	$P_{\alpha,k}=\e^{\lambda^\star E_k}-1$.
	Applying \eqref{eq:Mgr} at $t=\sigma_k$ and $t=\sigma_k+E_k$, we obtain	a.s.
	\[\e^{-\lambda^\star\sigma_k}(k+1)=W+\o(1),\qquad
	\e^{-\lambda^\star(\sigma_k+E_k)}(k+L_k+1)=W+\o(1).\]
	Since $W>0$ a.s., it readily follows that
	\begin{equation}
		\label{eq:atom_approx_0}
		\frac{L_k}{k+1}=P_{\alpha,k}+r_k(P_{\alpha,k}+1),
	\end{equation}
	where $r_k\to0$ a.s., and so
	\begin{equation}
		\label{eq:atom_approx}
		\frac{L_k}{n^{1+1/\alpha}}
		=\frac{(k+1)P_{\alpha,k}}{n^{1+1/\alpha}}+
		\frac{k+1}{n^{1+1/\alpha}}\,r_k(P_{\alpha,k}+1).
	\end{equation}
	The first term on the right-hand side is the second coordinate of the corresponding atom of
	$\widetilde\Lambda_n$. We now show that the second term
	is uniformly negligible on sets in $\Qc$.
	
	Fix $T>0$ as in \eqref{eq:Q}. For every $K\ge1$,
	\[\max_{k\le nT}\frac{k+1}{n^{1+1/\alpha}}
		|r_k|(P_{\alpha,k}+1)
		\le\frac{K}{n^{1+1/\alpha}}
		\max_{k<K}|r_k|(P_{\alpha,k}+1)\\
		+(T+1)\sup_{k\ge K}|r_k|\,
		\frac{\max_{l\le nT}(P_{\alpha,l}+1)}{n^{1/\alpha}}.\]
	For each fixed $K$, the first term tends to zero a.s.\ as $n\to\infty$.
	Since $P_{\alpha,l}+1$ is Pareto, the fraction in the second term
	is readily seen to be bounded in probability. Moreover,
	$\sup_{k\ge K}|r_k|\to0$ a.s.\ as $K\to\infty$.
	Hence, by letting first $n\to\infty$ and then $K\to\infty$,
	we see that the second term in \eqref{eq:atom_approx} vanishes
	in probability uniformly over $k\le nT$.
	Thus the difference between the second coordinates of corresponding
	atoms of $\widetilde\Lambda_n$ and $\Lambda_n$ is
	negligible on $B\in\Qc$, and $\widetilde\Lambda_n\xrightarrow{\mathrm{vd}}\Lambda$
	carries over to $\Lambda_n$.
\end{proof}
	
\begin{proof}[Proof of Theorem~\ref{thm:max_argmax}]
	Let
	\[\Pc_t=\supp\Lambda\cap\bigl([0,t]\times(0,\infty)\bigr).\]
	Define the argmax-max process $(X_\alpha,Y_\alpha)$ by
	\begin{equation}
	\label{eq:XY}
	\bigl(X_\alpha(t),Y_\alpha(t)\bigr)=
	\sum_{(x,y)\in\Pc_t}(x,y)\1\bigl\{y\ge y'\text{ for all }(x',y')\in\Pc_t\bigr\},\qquad t>0,
	\end{equation}
	with $X_\alpha(0)=Y_\alpha(0)=0$. The maximizing point is a.s.~unique since the intensity of $\Lambda$ is absolutely continuous. Hence, exactly one summand above is nonzero.
	
	We first prove that
	\begin{equation}
		\label{eq:ITtoXY}
		\left(\frac{I_{\lfloor n\cdot\rfloor}}n,\frac{T_{\lfloor n\cdot\rfloor}}{n^{1+1/\alpha}}\right)
		\xrightarrow{\mathrm d}(X_\alpha,Y_\alpha)
	\end{equation}
	in the $J_1$ topology on $D\bigl([0,\infty),\R^2\bigr)$. By
	the definition of $\Lambda_n$ in \eqref{eq:Lambda_n}, we have
	\[\left(\frac{I_{\lfloor nt\rfloor}}n,\frac{T_{\lfloor nt\rfloor}}{n^{1+1/\alpha}}	\right)
	=\bigl(X_{n,\alpha}(t),Y_{n,\alpha}(t)\bigr),\qquad t\ge0,\]
	where $(X_{n,\alpha},Y_{n,\alpha})$ is the argmax-max process associated
	with $\Lambda_n$, defined by the analogue of \eqref{eq:XY}, except that
	the indicator is supplemented by the condition $x\le u$ for every
	$(u,v)$ with $v=y$, thereby selecting the smallest argmax in the case
	of a tie.
	
	Applying the Skorokhod coupling to the convergence in
	Theorem~\ref{thm:point_process_conv}, we may assume that
	$\Lambda_n\to\Lambda$ vaguely a.s. We may also regard the second
	coordinate as taking values in $(0,\infty]$ and extend $\Lambda_n$ and
	$\Lambda$ by zero to $[0,\infty)\times\{\infty\}$. In view of the
	localization used in Theorem~\ref{thm:point_process_conv}, this extension
	does not affect the convergence. Hence, for fixed $\Delta>\delta>0$,
	\[[0,\Delta]\times\bigl[\tfrac12 Y_\alpha(\delta),\infty\bigr]\]
	is a.s.\ a compact $\Lambda$-continuity set. Proposition~3.13 in
	\cite{R08} therefore implies that the atoms of $\Lambda_n$ in this set
	converge a.s.~to the corresponding atoms of $\Lambda$.
	
	Since a.s.~no two atoms of $\Lambda$ have the same first or second
	coordinate, the record structure eventually agrees, and a piecewise
	linear time change matching the corresponding record times yields
	$(X_{n,\alpha},Y_{n,\alpha})\to(X_\alpha,Y_\alpha)$
	on $D([\delta,\Delta],\R^2)$. Moreover,
	$X_{n,\alpha}(t),X_\alpha(t)\le t$, while
	\[\sup_{0\le t\le\delta}Y_{n,\alpha}(t)
	=Y_{n,\alpha}(\delta)\to Y_\alpha(\delta),\]
	and $Y_\alpha(\delta)\to0$ as $\delta\downarrow0$. Thus the convergence
	extends to $D([0,\Delta],\R^2)$. Since $\Delta>0$ is arbitrary,
	\eqref{eq:ITtoXY} follows by analogy with Theorem 16.7 in \cite{B99}.
	
	It remains to prove that the argmax-max process $(X_\alpha,Y_\alpha)$
	defined by \eqref{eq:XY} is indeed the beta-decorated
	$\alpha$-Fr\'echet extremal process specified by the recursions
	\eqref{eq:Y_rec} and \eqref{eq:X_dist}. The recursive structure itself,
	apart from the asserted distributions, follows directly from
	\eqref{eq:XY}. Indeed, for $0\le s<t$, let
	$(B_{\alpha,s,t},M_{\alpha,s,t})$ be the highest atom of $\Lambda$ in
	$(s,t]\times(0,\infty)$. Then
	\[Y_\alpha(t)=\max\{Y_\alpha(s),M_{\alpha,s,t}\},\]
	and $X_\alpha(t)=X_\alpha(s)$ if the maximum does not change, whereas
	$X_\alpha(t)=B_{\alpha,s,t}$ otherwise. Moreover, by complete
	randomness of $\Lambda$,
	$(B_{\alpha,s,t},M_{\alpha,s,t})$ is independent of the
	restriction of $\Lambda$ to $[0,s]\times(0,\infty)$.
	
	It remains to identify the joint distribution of
	$(B_{\alpha,s,t},M_{\alpha,s,t})$. For $y>0$,
	\[\P\{\Lambda((s,t]\times(y,\infty))=0\}=\exp\biggl\{
	-\int_s^t\int_y^\infty\alpha x^\alpha z^{-\alpha-1}\,\dd x\,\dd z
	\biggr\}=\e^{-\|g_{s,t}\|_\alpha^\alpha y^{-\alpha}}.\]
	Hence, by the Mecke formula (see, e.g., Theorem~4.1 in
	\cite{LP18}), for every non-negative Borel function $h$,
	\[\begin{aligned}
		\E h(B_{\alpha,s,t},M_{\alpha,s,t})
		&=\E\sum_{(x,y)\in\supp\Lambda}
		h(x,y)\1\{s<x\le t\}\1\{\Lambda((s,t]\times(y,\infty))=0\}\\
		&=\int_s^t\int_0^\infty
		h(x,y)\e^{-\|g_{s,t}\|_\alpha^\alpha y^{-\alpha}}
		\alpha x^\alpha y^{-\alpha-1}\,\dd x\,\dd y.
	\end{aligned}\]
	Thus, $(B_{\alpha,s,t},M_{\alpha,s,t})$ has joint density
	\[\alpha x^\alpha y^{-\alpha-1}\e^{-\|g_{s,t}\|_\alpha^\alpha y^{-\alpha}}\1_{(s,t]}(x),\]
	which is easily seen to factorize into the density of $B_{\alpha,s,t}$ and that of
	$\|g_{s,t}\|_\alpha F_{\alpha,s,t}$. This factorization also shows that the two components are independent;
	as noted above, the same can alternatively be deduced from Theorem~8.2	in \cite{MS16}.
	\end{proof}
	
\section{\texorpdfstring{$\ell_k$}{ell\_k}-long-lasting leaves}
\label{sec:5}

Under \eqref{eq:M} and \eqref{eq:L}, Proposition~\ref{prop:L} shows that the lifetime of vertex $n$ as a leaf is typically of order $n$. At the same time, Theorem~\ref{thm:max_argmax} and
Corollary~\ref{cor:IT} show that the maximum leaf lifetime among the
vertices $0,\ldots,n$ is of order $n^{1+1/\alpha}$, where
$\alpha=f(0)/\lambda^\star$. In this section, we analyze the asymptotic
behavior of the number of vertices that remain leaves for unusually long
periods of time.

Fix a deterministic positive threshold sequence
$(\ell_k,k\ge1)$. We call vertex $k\ge1$ an $\ell_k$-long-lasting leaf
if $L_k\ge\ell_k$; since $L_0=1$, vertex $0$ is excluded from consideration in what follows. We also introduce the random set of long-lasting leaves
\[\L=\{k\ge1:L_k\ge\ell_k\}\]
and start with a criterion for the a.s.\ finiteness of $\L$.

\begin{proposition}
	\label{prop:cardL}
	Under \eqref{eq:M} and \eqref{eq:L}, the set $\L$ is a.s.\
	finite if $\sum_{k=1}^\infty(k/\ell_k)^\alpha<\infty$,
	and a.s.\ infinite otherwise.
\end{proposition}

For the power threshold $\ell_k=k^\beta$, this yields
$\card\L=\infty$ a.s.\ if and only if
$\beta\le1+1/\alpha$, which agrees with the normalization of $T_n$ in
Theorem~\ref{thm:max_argmax} and Corollary~\ref{cor:IT}.

\begin{proof}[Proof of Proposition \ref{prop:cardL}]
	As in the proof of Theorem~\ref{thm:point_process_conv}, let
	$P_{\alpha,k}=\e^{\lambda^\star E_k}-1$. Then, by \eqref{eq:atom_approx_0},
	$1+\frac{L_k}{k+1}=(1+P_{\alpha,k})(1+r_k)$,
	where $r_k\to0$ a.s. Hence, a.s.\ for all large $k$,
	\begin{equation}
		\label{eq:L-P-comparison}
		\frac12(1+P_{\alpha,k})
		\le 1+\frac{L_k}{k+1}\le
		\frac32(1+P_{\alpha,k}).
	\end{equation}
	
	Note that the series
	$\sum_{k=1}^\infty\bigl(\frac k{\ell_k}\bigr)^\alpha$
	and
	$\sum_{k=1}^\infty
	\bigl(1+\frac{\ell_k}{k+1}\bigr)^{-\alpha}$
	converge or diverge together. Moreover, by \eqref{eq:L-P-comparison},
	a.s.\ for large $k$,
	\[1+P_{\alpha,k}\ge
	2\biggl(1+\frac{\ell_k}{k+1}\biggr)
	\quad\Longrightarrow\quad L_k\ge\ell_k
	\quad\Longrightarrow\quad
	1+P_{\alpha,k}\ge\frac23\biggl(1+\frac{\ell_k}{k+1}\biggr).\]
	The probabilities of the leftmost events are
	$2^{-\alpha}\bigl(1+\frac{\ell_k}{k+1}\bigr)^{-\alpha}$,
	while, in the convergent case and for large $k$, the probabilities
	of the rightmost events are
	$\bigl(\frac23\bigr)^{-\alpha}
	\bigl(1+\frac{\ell_k}{k+1}\bigr)^{-\alpha}$.
	Since the $P_{\alpha,k}$ are independent, the claim follows from the
	Borel--Cantelli lemma.
\end{proof}

We now study Gaussian fluctuations at the critical power threshold,
that is, at the boundary between finitely and infinitely many
long-lasting leaves. Here \eqref{eq:M} and \eqref{eq:L} alone turn out
to be insufficient, since we will need a rate of convergence in the
Malthusian growth theorem \eqref{eq:Mgr}. Fortunately, a bound
from \cite{BBC23} provides the required control for the broad
class of attachment functions introduced in
Proposition~\ref{prop:cond}.

\begin{theorem}
	\label{thm:FCLT}
	Under the assumptions of Proposition~\ref{prop:cond}, let
	$\ell_k=c k^{1+1/\alpha}$ for some $c>0$ and set
	$S_n=\card(\L\cap\{1,\ldots,n\})$. Then, as $n\to\infty$,
	\begin{equation}
		\label{eq:FCLT}
		\biggl(\frac{S_{\lfloor n^t\rfloor}-c^{-\alpha}t\log n}
		{\sqrt{c^{-\alpha}\log n}}\biggr)_{t\ge0}\!\!\!
		\xrightarrow{\mathrm d}(B(t))_{t\ge0}
	\end{equation}
	in the $J_1$ topology on $D([0,\infty),\R)$, where $B$ is
	standard Brownian motion.
\end{theorem}

Arguments similar to those in the proof below can be applied to a
broader class of threshold sequences close to the critical power.
We do not pursue this extension, since the centering and normalization
would become considerably less explicit, whereas the critical threshold
is, in our view, of particular interest in its own right. As for lower power thresholds, a more substantial extension would require stronger assumptions on the attachment function $f$, at the expense of the universality of our result.

\begin{proof}[Proof of Theorem~\ref{thm:FCLT}]
	We use the notation $\sigma_k$, $E_k$, $P_{\alpha,k}$, and $N$
	introduced above. First note that, under the assumptions of
	Proposition~\ref{prop:cond}, the $\o(1)$ in \eqref{eq:Mgr}
	can be replaced by $\O(t^{-2})$; in fact, even an exponential
	rate is available, but this will suffice for us. Indeed,
	Theorem~3.5 in \cite{BBC23}, applied to the population-counting
	characteristic $\varphi=\1_{[0,\infty)}$, yields exponential decay of
	$\E|\e^{-\lambda^\star t}N(t)-W|$ as $t\to\infty$.
	By taking $t_i=i^{1/3}$, a standard application of
	Markov's inequality and the Borel--Cantelli lemma shows that
	$\e^{-\lambda^\star t_i}N(t_i)-W$ itself decays exponentially a.s.
	Now for $t_i\le t\le t_{i+1}$, the monotonicity of $N$ yields
	\begin{align*}
	\e^{-\lambda^\star t_i}N(t_i)\cdot\e^{\lambda^\star(t_i-t_{i+1})}
	&=\e^{-\lambda^\star t_{i+1}}N(t_i)\le\e^{-\lambda^\star t}N(t)
	\\&\le\e^{-\lambda^\star t_i}N(t_{i+1})=
	\e^{-\lambda^\star t_{i+1}}N(t_{i+1})\cdot\e^{\lambda^\star(t_{i+1}-t_i)},
	\end{align*}
	and hence, since $t_{i+1}-t_i$ is of order $t_i^{-2}$,
	\begin{equation*}
		\e^{-\lambda^\star t}N(t)=W+\O(t^{-2})\quad\text{a.s.}
	\end{equation*}
	
	As in the proof of Proposition~\ref{prop:L}, applying this
	at $t=\sigma_k$ and $t=\sigma_k+E_k$, and using
	$\sigma_k\asymp\log k$ from \eqref{eq:siglog} and $W>0$, gives
	\eqref{eq:atom_approx_0} with $r_k=\O((\log k)^{-2})$ a.s.
	In other words,
	\begin{equation}
	\label{eq:logs}
	\biggl|\log\biggl(1+\frac{L_k}{k+1}\biggr)
	-\log(1+P_{\alpha,k})\biggr|=\O\bigl((\log k)^{-2}\bigr)
	\quad\text{a.s.}
	\end{equation}
	Since the constant implicit in $\O$ may be random, we shall henceforth
	use \eqref{eq:logs} with $=\O\bigl((\log k)^{-2}\bigr)$ replaced by
	$\le(\log k)^{-3/2}$ a.s.\ for large $k$. Now set 
	\begin{gather*}
	\I_k=\1\{L_k\ge\ell_k\}=\1\biggl\{\log\biggl(1+\frac{L_k}{k+1}\biggr)
	\ge\log\biggl(1+\frac{\ell_k}{k+1}\biggr)\biggr\},\\
	\widetilde\I_k=\1\{(k+1)P_{\alpha,k}\ge\ell_k\}=\1\biggl\{\log(1+P_{\alpha,k})\ge
	\log\biggl(1+\frac{\ell_k}{k+1}\biggr)\biggr\}
	\end{gather*}
	and show that these indicators agree eventually a.s.	
	Indeed, if $\I_k\ne\widetilde\I_k$, then
	$\log(1+\ell_k/(k+1))$ lies between the two logarithms on the left-hand sides of the inequalities. Hence, by \eqref{eq:logs},
	\[\I_k\ne\widetilde\I_k
	\quad\Longrightarrow\quad
	\biggl|\log(1+P_{\alpha,k})
	-\log\biggl(1+\frac{\ell_k}{k+1}\biggr)\biggr|
	\le(\log k)^{-3/2}\]
	a.s.\ for large $k$. Since $\log(1+P_{\alpha,k})\sim\mathsf{Exp}(\alpha)$,
	the probability of the event on the right equals
	\[\biggl(1+\frac{\ell_k}{k+1}\biggr)^{-\alpha}
	\bigl(\e^{\alpha(\log k)^{-3/2}}
	-\e^{-\alpha(\log k)^{-3/2}}\bigr)
	=\O\bigl(k^{-1}(\log k)^{-3/2}\bigr).\]
	These probabilities are summable, so the Borel--Cantelli lemma
	yields $\I_k=\widetilde\I_k$ eventually a.s. Thus, it suffices to prove
	\eqref{eq:FCLT} with $S_{\lfloor n^t\rfloor}=\sum_{k=1}^{\lfloor n^t\rfloor}\I_k$ replaced by
	$\widetilde S_{\lfloor n^t\rfloor}
	=\sum_{k=1}^{\lfloor n^t\rfloor}\widetilde\I_k$.
	
	Unlike $\I_k$, the indicators $\widetilde\I_k$ are independent.
	Furthermore,
	\[\E\widetilde\I_k=\biggl(1+\frac{\ell_k}{k+1}\biggr)^{-\alpha},
	\quad
	\var\widetilde\I_k=
	\biggl(1+\frac{\ell_k}{k+1}\biggr)^{-\alpha}\!\!-
	\biggl(1+\frac{\ell_k}{k+1}\biggr)^{-2\alpha},\]
	and both are asymptotic to $c^{-\alpha}k^{-1}$ as $k\to\infty$. Moreover, it is straightforward to check that
	$\sum_{k=1}^\infty|\E\widetilde\I_k-c^{-\alpha}k^{-1}|<\infty$ and $\sum_{k=1}^\infty|\var\widetilde\I_k-c^{-\alpha}k^{-1}|<\infty$.
	Thus,
	\[\sum_{k=1}^{\lfloor n^t\rfloor}\E\widetilde\I_k
	=c^{-\alpha}t\log n+\O(1),\qquad
	\sum_{k=1}^{\lfloor n^t\rfloor}\var\widetilde\I_k
	=c^{-\alpha}t\log n+\O(1),\]
	uniformly for $t\ge0$. Since $\widetilde\I_k$ are bounded, the Lindeberg condition follows automatically, and the functional central limit theorem for independent
	non-identically distributed random variables yields the desired convergence.
\end{proof}

\section*{Acknowledgments}

Part of this research was carried out while the second author
was affiliated with the Institute of Mathematical Statistics
and Actuarial Science at the University of Bern,
Switzerland. He is deeply grateful to his former colleagues
there, whose kindness and support made it such a welcoming
and stimulating place to work.

\section*{A note on the use of AI}

The idea behind this paper originated quite some time ago
in an attempt to relate certain aspects of the evolution
of random recursive trees to random point measures.
Active work began about a year and a half ago
and continued intermittently, initially without AI assistance.

During the final stages, we used ChatGPT (OpenAI, GPT-6).
This collaboration led to several improve\-ments in the paper:
\begin{itemize}
	\item ChatGPT drew our attention to
	Theorem~6.11 in \cite{HJ17}, which allowed us
	to simplify the proof of Proposition~\ref{prop:EL}.
	
	\item Through extensive discussions with ChatGPT,
	we developed Proposition~\ref{prop:cond},
	which provides explicit conditions covering a broad
	class of attachment functions and extends earlier
	results from \cite{GVCH17} and \cite{RTV07}.
	
	\item ChatGPT also pointed out Theorem~3.5 in
	\cite{BBC23}, which, through further joint exploration,
	enabled us to substantially strengthen the
	result now stated as Theorem~\ref{thm:FCLT}.
	An earlier version had been proved
	as a one-dimensional limit theorem for affine
	preferential-attachment trees, using a completely
	different approach based on Stein's method.
\end{itemize}

ChatGPT was also used for language editing.
All mathematical statements, proofs, and references were
carefully and independently verified by the authors,
who take full responsibility for the content of the paper.

\authorinfo{Oleksii Galganov}
{Department of Mathematical Analysis and Probability Theory,
Igor Sikorsky Kyiv Polytechnic Institute, Prospect Beresteiskyi~37, 03056 Kyiv, Ukraine}{galganov.oleksii@edu.kpi.ua}
\bigskip

\authorinfo{Andrii Ilienko}
{Faculty of Pure and Applied Mathematics,
Wroc{\l}aw University of Science and Technology,
Wybrze\.ze Wyspia\'nskiego 27,
50-370 Wroc{\l}aw, Poland}
{andrii.ilienko@pwr.edu.pl}
\end{document}